\documentclass{article}

\usepackage{subfigure}
\usepackage{color}
\usepackage{amsmath,amsfonts}
\usepackage{verbatim}
\usepackage{latexsym}
\usepackage[utf8]{inputenc} 
\usepackage{amssymb}
\usepackage{mathrsfs}
\usepackage{graphicx}
\usepackage{threeparttable}

\usepackage[noblocks]{authblk}

\usepackage{caption}
\usepackage{epstopdf}

\date{ }

\begin{document}
	\newcommand{\E}{\mathbb{E}}
	\newcommand{\PP}{\mathbb{P}}
	\newcommand{\RR}{\mathbb{R}}
	
	\newcommand{\Dt}{\D t}
	\newcommand{\bX}{\bar X}
	\newcommand{\bx}{\bar x}
	\newcommand{\by}{\bar y}
	\newcommand{\bp}{\bar p}
	\newcommand{\bq}{\bar q}

	\newtheorem{thm}{Theorem}[section]
	\newtheorem{lemma}[thm]{Lemma}
	\newtheorem{coro}[thm]{Corollary}
	\newtheorem{defn}[thm]{Definition}
	\newtheorem{assp}[thm]{Assumption}
	\newtheorem{expl}[thm]{Example}
	\newtheorem{prop}[thm]{Proposition}
	\newtheorem{rmk}[thm]{Remark}



\title{Multi-level Monte Carlo methods with the truncated Euler-Maruyama scheme for stochastic differential equations}
\author{Qian Guo\textsuperscript{a}
Wei Liu\textsuperscript{a}
Xuerong Mao\textsuperscript{b}
and Weijun Zhan\textsuperscript{a}
\thanks{Corresponding author. Email: weijunzhan@hotmail.com\\
\ \ \textsuperscript{a}Department of Mathematics, Shanghai Normal University, Shanghai, China;\\
\ \ \textsuperscript{b}Department of Mathematics and Statistics, University of Strathclyde, Glasgow, UK}
}
\maketitle

\begin{abstract}
The truncated Euler-Maruyama method is employed together with the Multi-level Monte Carlo method to approximate expectations of some functions of solutions to stochastic differential equations (SDEs). The convergence rate and the computational cost of the approximations are proved, when the coefficients of SDEs satisfy the local Lipschitz and Khasminskii-type conditions. Numerical examples are provided to demonstrate the theoretical results.\\

\end{abstract}
%
%

\section{Introduction}

Stochastic differential equations (SDEs) have been broadly discussed and applied as a powerful tool to capture the uncertain phenomenon in the evolution of systems in many areas \cite{Allen2007,Gardiner2004,Mao2008,Oksendal2003,PH2006}. However, the explicit solutions of SDEs can rarely be found. Therefore, the numerical approximation  becomes an essential approach in the applications of SDEs. Monographs \cite{KP1992,MT2004} provide detailed introductions and discussions to various classic methods.
\par
Since the non-linear coefficients have been widely adapted in SDE models \cite{Ait1996,GGHMP2011,NBC2016}, explicit numerical methods that have good convergence property for SDEs with non-global Lipschitz drift and diffusion coefficients are of interest to many researchers and required by practitioners. The authors in \cite{HMS2002} developed a quite general approach to prove the strong convergence of numerical methods for nonlinear SDEs. The approach to prove the global strong convergence via the local convergence for SDEs with non-global Lipschitz coefficients was studied in \cite{TZ2013}. More recently, the taming technique was developed to handle the non-global Lipschitz coefficients \cite{HJ2015,HJK2012}. Simplified proof of the tamed Euler method and the tamed Milstein method can be found in \cite{Sabanis2013} and \cite{WG2013}, respectively. The truncated Euler-Maruyama (EM) method was developed in \cite{Mao2015,Mao2016}, which is also targeting on SDEs with non-global Lipschitz coefficients. Explicit methods for nonlinear SDEs that preserve positivity can be found in, for example \cite{Halidias2015,LM2013}. {\color{blue} A modified truncated EM  method that preserves the asymptotic stability and boundedness of the nonlinear SDEs was presented in \cite{GLMY2017}.}
\par
Compared to the explicit methods mentioned above, the methods with implicit term have better convergence property in approximating non-global Lipschitz SDEs with the trade-off of the relatively expensive computational cost. We just mention a few of the works \cite{Hu1996,SMHP2011,YHJ2014} and the references therein.
\par
In many situations, the expected values of some functions of the solutions to SDEs are also of interest. To estimate the expected values, the classic Monte-Carlo method is a good and natural candidate. More recently, Giles in \cite{Giles2008b,Giles2008a} developed the Multi-level Monte Carlo (MLMC) method, which improves the convergence rate and reduces the computational cost of estimating expected values. {\color{blue} A detailed survey of recent developments and applications of the MLMC method can be found in \cite{Giles2015}. To complement \cite{Giles2015}, we only mention some new developments that are not included in \cite{Giles2015}. Under the global Lipschitz and linear growth conditions, the MLMC method combined with the EM method applied to SDEs with small noise is often found to be the most efficient option \cite{Anderson2016}. The MLMC method with the adaptive EM method was designed for solving SDEs driven by L\'{e}vy process \cite{Dereich2016a,Dereich2016b}. The MLMC method was applied to SDEs driven by Poisson random measures by means of coupling with the split-step implicit tau-leap at levels. However,} the classic EM method with the MLMC method has been proved divergence to SDEs with non-global Lipschitz coefficients \cite{HJK2013}. So it is interesting to investigate the combinations of the MLMC method with those numerical methods developed particularly for SDEs with non-global Lipschitz coefficients. In \cite{HJK2013}, the tamed Euler method was combined with the MLMC method to approximate expectations of some nonlinear functions of solutions to some nonlinear SDEs.
\par
In this paper, we embed the MLMC method with the truncated EM method and study the convergence and the computational cost of this combination to approximate expectations of some nonlinear functions of  solutions to SDEs with non-global Lipschitz coefficients.
\par
In \cite{Mao2016}, the truncated EM method has been proved to converge to the true solution with the order $1/2-\varepsilon$ for any arbitrarily small $\varepsilon > 0$. The plan of this paper is as follows. Firstly, we make some modifications of  Theorem 3.1 in \cite{Giles2008a} such that the modified theorem is able to cover the truncated EM method. Then, we use the modified theorem to prove the convergence and the computational cost of the MLMC method with the truncated EM method. At last, numerical examples for SDEs with non-global Lipschitz coefficients and expectations of nonlinear functions are given to demonstrate the theoretical results.
\par
This paper is constructed as follows. Notations, assumptions and some existing results about the truncated EM method and the MLMC method are presented in Section 2. Section 3 contains the main result on the computational complexity. A numerical example is provided in Section 4 to illustrate theoretical results.
In the appendix, we give the proof of the theorem in Section 3.

\section{Mathematical Preliminary}

Throughout this paper, unless otherwise specified, we let $(\mathit{\Omega},\mathscr{F},\mathbb{P})$ be a complete probability
space with a filtration $\{\mathscr{F}_t\}_{t\geq 0}$ satisfying the usual condition (that is, it is right continuous and increasing
while $\mathscr{F}_0$ contains all $\mathbb{P}-$null sets). Let $\mathbb{E}$ denote the expectation corresponding to $\mathbb{P}$.
Let $B(t)$ be an $m$-dimensional Brownian motion defined on the space. If $A$ is a vector or matrix, its transpose is denoted by $A^T$.
If $x\in \mathbb{R}^d$, then $|x|$ is the Euclidean norm. If $A$ is a matrix, we let $|A|=\sqrt{\mathrm{trace}(A^TA)}$ be its trace norm.
If $A$ is a symmetric matrix, denote by $\lambda_{\max}(A)$ and $\lambda_{\min}(A)$ its largest and smallest eigenvalue, respectively.
Moreover, for two real numbers $a$ and $b$, set $a\vee b=\max(a,b)$ and $a\wedge b=\min(a,b)$. If $G$ is a set, {\color{blue} {its indicator function
is denoted}} by $I_G(x)=1$ if $x\in G$ and $0$ otherwise.\par
Here we consider an SDE
\begin{equation}\label{underlyingSDE}
	dX(t)=\mu(X(t))dt+\sigma(X(t))dB(t)
\end{equation}
on $t\geq 0$ with the initial value $X(0)=X_0\in \mathbb{R}^d$, where
{\color{blue}
\begin{equation*}
	\mu : \mathbb{R}^d \rightarrow \mathbb{R}^d \ \ \ \text{and}\ \ \ \sigma : \mathbb{R}^d \rightarrow \mathbb{R}^{d\times m}.
\end{equation*}\par
}
When the coefficients obey the global Lipschitz condition, the strong convergence of numerical methods for SDEs has been well studied \cite{KP1992}.
When the coefficients $\mu$ and $\sigma$ are locally Lipschitz continuous without the linear growth condition,
Mao \cite{Mao2015,Mao2016} recently developed the truncated EM method.
To make this paper self-contained, we give a brief review of this method firstly.
\par
We first choose a strictly increasing continuous function {\color{blue} $\omega : \mathbb{R}_+\rightarrow \mathbb{R}_+$ such that $\omega(r)\rightarrow  \infty$ as $r\rightarrow \infty$} and
\begin{equation}\label{assp3}
	\sup\limits_{|x|\leq u}(|\mu(x)|\vee|\sigma(x)|)\leq \omega(u),\ \ \ \forall u\geq 1.
\end{equation}
Denote by $\omega^{-1}$ the inverse function of $\omega$ and we see that $\omega^{-1}$ is a strictly increasing continuous function from
$[\omega(0),\infty)$ to $R_+$.
We also choose a number $s_l^*\in (0,1]$ and a strictly decreasing function {\color{blue} $h : (0,s_l^*]\rightarrow (0,\infty)$}
such that
\begin{equation}\label{assp4}
	h(s_l^*)\geq \omega(2),\ \ \ \lim\limits_{s_l\rightarrow 0}h(s_l)=\infty \ \  \text{and} \ \ s_l^{1/4}h(s_l)\leq1,\ \ \ \forall s_l\in(0,s_l^*].
\end{equation}
For a given stepsize $s_l\in(0,1)$, let us define the truncated functions
\begin{equation*}
\mu_{s_l}(x)=\mu\left((|x|\wedge \omega^{-1}(h(s_l)))\frac{x}{|x|}\right)\ \ \textup{and}\ \
\ \sigma_{s_l}(x)=\sigma\left((|x|\wedge \omega^{-1}(h(s_l)))\frac{x}{|x|}\right)
\end{equation*}
for $x\in\mathbb{R}^d$, where we set $x/|x|=0$ when $x=0$.
Moreover, let $\overline{X}_{s_l}(t)$ denote the approximation to $X(t)$ using the truncated EM method with time step size $s_l=M^{-l}T$ for  $l=0,1,\ldots,L$.
The numerical solutions $\overline{X}_{s_l}(t_k)$ for $t_k=ks_l$ are formed by setting
$\overline{X}_{s_l}(0)=X_0$ and computing
\begin{equation}\label{truncated EM}
\overline{X}_{s_l}(t_{k+1})=\overline{X}_{s_l}(t_{k})+\mu_{s_l}(\overline{X}_{s_l}(t_{k}))s_l+\sigma_{s_l}(\overline{X}_{s_l}(t_{k}))\Delta B_k
\end{equation}
for $k=0,1,\ldots,$ where $\Delta B_k=B(t_{k+1})-B(t_k)$ is the Brownian motion increment.

Now we give some assumptions to guarantee that the truncated EM solution (\ref{truncated EM}) will converge to the true solution to the SDE (\ref{underlyingSDE}) in the strong sense. \par
{\color{blue}
\begin{assp}\label{assp1}
 The coefficients $\mu$ and $\sigma$ satisfy the local Lipschitz condition that for any real number $R>0$,there exists a $K_R>0$ such that
\begin{equation}
	|\mu(x)-\mu(y)|\vee |\sigma(x)-\sigma(y)|\leq K_R|x-y|
\end{equation}
for all $x,y\in R^d$ with $|x|\vee|y|\leq R$.
\end{assp}
}
\begin{assp}\label{assp2}
 The coefficients  $\mu$ and $\sigma$ satisfy the Khasminskii-type condition that there exists a pair of constants $p>2$ and $K>0$ such that
\begin{equation}
	x^T\mu(x)+\frac{p-1}{2}|\sigma(x)|^2\leq K(1+|x|^2)
\end{equation}
for all $x\in \mathbb{R}^d$.
\end{assp}

\begin{assp}\label{assp5}
 There exists a pair of constants $q\geq2$ and $H_1>0$ such that
\begin{equation}
	(x-y)^T(\mu(x)-\mu(y))+\frac{q-1}{2}|\sigma(x)-\sigma(y)|^2\leq H_1|x-y|^2
\end{equation}
for all $x,y\in \mathbb{R}^d$.
\end{assp}
\begin{assp}\label{assp6}
There exists a pair of positive constants $\rho$ and $H_2$ such that
\begin{equation}
	|\mu(x)-\mu(y)|^2\vee|\sigma(x)-\sigma(y)|^2\leq H_2(1+|x|^\rho+|y|^\rho)|x-y|^2
\end{equation}
for all $x,y\in \mathbb{R}^d$.
\end{assp}
{\color{blue} Let $f(X(t))$ denote a payoff function of the solution to some SDE driven by a given Brownian path $B(t)$. In this paper, we need $f$ satisfies the following assumption.}
\begin{assp}
There exists a constant $c>0$ such that
\begin{equation}\label{poly grown condition}
	|f(x)-f(y)|\leq c(1+|x|^c+|y|^c)|x-y|
\end{equation}
for all $x,y\in \mathbb{R}^d$.

\end{assp}

Using the idea in \cite{Giles2008b,Giles2008a}, the expected value of $f(\overline{X}_{s_l}(t))$
can be decomposed in the following way
\begin{equation}
\label{decomp}
	\mathbb{E}[f(\overline{X}_{s_L}(T))]=\mathbb{E}[f(\overline{X}_{s_0}(T))]+\sum\limits_{l=1}^L \mathbb{E}[f(\overline{X}_{s_l}(T))-f(\overline{X}_{s_{l-1}}(T))].
\end{equation}\par
Let $Y_0$ be an estimator for $\mathbb{E}[f(\overline{X}_{s_0}(T))]$ using $N_0$ samples. Let $Y_l$ be an estimator for $\mathbb{E}[f(\overline{X}_{s_l}(T))-f(\overline{X}_{s_{l-1}}(T))]$ using $N_l$ paths such that
$$Y_l=\frac{1}{N_l}\sum_{i=1}^{N_l}\left[f(\overline{X}_{s_l}^{(i)}(T))-f(\overline{X}_{s_{l-1}}^{(i)}(T))\right].$$
The multi-level method independently estimates each
of the expectations on the right-hand side of (\ref{decomp}) such that the computational complexity can be minimized, see \cite{Giles2008a} for more details. \par

\section{Main Results}

In this section, {Theorem 3.1} in \cite{Giles2008a} is slightly generalised. Then the convergence rate and computational complexity of the truncated EM method combined with the MLMC method are studied.

\subsection{Generalised theorem for the MLMC method}
\begin{thm}\label{generalthmMLMC}
If there exist independent estimators $Y_l$ based on $N_l$ Monte Carlo samples, and positive constants $\alpha,\beta,c_1,c_2,c_3$ such that \par
\begin{enumerate}
\item $ \mathbb{E}[f(\overline{X}_{s_l}(T))-f(X(T))]\leq c_1s_l^\alpha$,\par
\item \begin{eqnarray*}\mathbb{E}[Y_l]=
		\begin{cases}
			\mathbb{E}[f(\overline{X}_{s_0}(T))], &l=0,\cr \mathbb{E}[f(\overline{X}_{s_l}(T))-f(\overline{X}_{s_{l-1}}(T))], & l>0,
		\end{cases}
	\end{eqnarray*}
\item $Var[Y_l]\leq c_2N_l^{-1}s_l^\beta$,
\item the computational complexity of $Y_l$, denoted by $C_l$, is bounded by
$$C_l\leq c_3N_ls_l^{-1},$$
\end{enumerate}
then there exists a positive constant $c_4$ such that for any $\varepsilon<e^{-1}$ the multi-level estimator
$$Y=\sum\limits_{l=0}^{L}Y_l $$
has a mean square error (MSE)
$$MSE\equiv \mathbb{E}\Big[\Big(Y-\mathbb{E}[f(X(T))]\Big)^2\Big]<\varepsilon^2.$$

Furthermore, the upper bound of computational complexity of $Y$, denoted by $C$, is given by
\begin{eqnarray*}C\leq
\begin{cases}
			c_3(2c_5^2c_2+\frac{M^2}{M-1}(\sqrt{2}c_1)^{1/\alpha}){\varepsilon^{-1/\alpha}}, &\alpha\leq (-log{\varepsilon})/{log[(log\varepsilon/\varepsilon)^2]},\cr
			c_3(2c_5^2c_2+\frac{M^2}{M-1}(\sqrt{2}c_1)^{1/\alpha}){\varepsilon^{-2}(log\varepsilon)^2}, &\alpha> (-log{\varepsilon})/{log[(log\varepsilon/\varepsilon)^2]}
\end{cases}
\end{eqnarray*}
for $\beta=1$,
\begin{eqnarray*}C\leq
		\begin{cases}
			c_3[2c_2T^{\beta-1}(1-M^{-(\beta-1)/2})^{-2}+\frac{M^2}{M-1}(\sqrt{2}c_1)^{1/\alpha}]\varepsilon^{-2}, &\alpha\geq \frac{1}{2},\cr
			c_3[2c_2T^{\beta-1}(1-M^{-(\beta-1)/2})^{-2}+\frac{M^2}{M-1}(\sqrt{2}c_1)^{1/\alpha}]\varepsilon^{-1/\alpha},& {\color{blue}\alpha< \frac{1}{2}}
		\end{cases}
\end{eqnarray*}
for $\beta>1$, and
\begin{eqnarray*}C\leq
		\begin{cases}
			c_3[2c_2(\sqrt{2}c_1)^{(1-\beta)/\alpha}M^{1-\beta}(1-M^{-(1-\beta)/2})^{-2}
			+\frac{M^2}{M-1}(\sqrt{2}c_1)^{1/\alpha}]\varepsilon^{-2-(1-\beta)/\alpha},&\beta\leq 2\alpha, \cr
			 c_3[2c_2(\sqrt{2}c_1)^{(1-\beta)/\alpha}M^{1-\beta}(1-M^{-(1-\beta)/2})^{-2}
			+\frac{M^2}{M-1}(\sqrt{2}c_1)^{1/\alpha}]\varepsilon^{-1/\alpha},&\beta> 2\alpha
		\end{cases}
	\end{eqnarray*}
 for $0<\beta<1$.
\end{thm}
The proof is in Appendix.
\begin{rmk}
The main difference of Theorem \ref{generalthmMLMC} and {Theorem 3.1} in \cite{Giles2008a} lies in the first condition. In \cite{Giles2008a}, one needs $\alpha\geq\frac{1}{2}$. In this paper, this requirement is weaken by any $\alpha>0$.
\end{rmk}

\subsection{Specific theorem for truncated Euler with the MLMC}
Next we consider the multi-level Monte Carlo path simulation with truncated
EM method and discuss their computational complexity using Theorem  {\ref{generalthmMLMC}}.

From Theorem 3.8 in \cite{Mao2016}, under {\color{blue} Assumptions {\ref{assp1}}-{\ref{assp6}}},
for every small $s_l\in(0,s_l^{\ast})$, where $s_l^{\ast}\in(0,1)$ and for any real number $T>0$, we have
\begin{equation}\label{strong convergence}
	\mathbb{E}|X(T)-\overline{X}_{s_l}(T)|^{\overline{q}}\leq c~s_l^{\overline{q}/2}~(h(s_l))^{\overline{q}},
\end{equation}
for $\overline{q}\geq 2$.
If $\overline{q}=1$, by using the Holder inequality, we also know that
\begin{equation*}
	\mathbb{E}|X(T)-\overline{X}_{s_l}(T)|\leq {(\mathbb{E}|X(T)-\overline{X}_{s_l}(T)|^2)}^{1/2}
	\leq (c s_l (h(s_l))^2)^{1/2}=c s_l^{1/2} h(s_l),
\end{equation*}
so we can obtain
\begin{equation}\label{weak convergence}
	\begin{aligned}
		&\mathbb{E}[|f(\overline{X}_{s_l}(T))-f(X(T))|]\\
		&\leq \mathbb{E}[c(1+|\overline{X}_{s_l}(T)|^c+|X(T)|^c)|\overline{X}_{s_l}(T)-X(T)|]
		\leq c (\mathbb{E}|\overline{X}_{s_l}(T)-X(T)|^2)^{1/2}\\
		&\leq c s_l^{1/2} h(s_l)
	\end{aligned}
\end{equation}
with the polynomial growth condition {(\ref{poly grown condition}}). {This implies that $\alpha=1/4$ for the truncated EM scheme.}
\par
Next we consider the variance of $Y_l$.
It follows that
\begin{equation}\label{var}
	Var[f(\overline{X}_{s_l}(T))-f(X(T))]\leq \mathbb{E}[\big(f(\overline{X}_{s_l}(T))-f(X(T))\big)^2]\leq c s_l (h(s_l))^2
\end{equation}
{\color{blue}
using (\ref{poly grown condition}) and (\ref{strong convergence}).
}
In addition, it can be noted that
\begin{equation*}
\begin{aligned}
& f(\overline{X}_{s_l}(T))-f(\overline{X}_{s_{l-1}}(T))=[f(\overline{X}_{s_l}(T))-f(X(T))]-[f(\overline{X}_{s_{l-1}}(T))-f(X(T))],
\end{aligned}
\end{equation*}
{
thus we have
\begin{equation*}
\begin{aligned}
&Var[f(\overline{X}_{s_l}(T))-f(\overline{X}_{s_{l-1}}(T))]\\
{\color{blue} \leq } &{\color{blue} \left( \sqrt{Var[f(\overline{X}_{s_l}(T))-f(X(T))]}+\sqrt{Var[f(\overline{X}_{s_{l-1}}(T))-f(X(T))]}\right)^2}\\
\leq & c s_l (h(s_l))^2+c s_{l-1} (h(s_{l-1}))^2\\
\leq & c s_l^{1/2},
\end{aligned}
\end{equation*}
where the fact $s_l^{\frac{1}{4}}~h(s_l)\leq 1$  from (\ref{assp4}) is used.
{\color{blue} Now we have $$Var[Y_l]=N_{l}^{-1}Var\left[f(\overline{X}_{s_l}^{(i)}(T))-f(\overline{X}_{s_{l-1}}^{(i)}(T))\right]\leq c N_l^{-1} s_l^{1/2}.$$}
So we have $\beta=1/2$ for the truncated EM method.

 }
{
According to the Theorem 3.1, it is easy to find that the upper bound of the computational complexity of $Y$ is
$$\left[4 c_1^2 c_2 c_3\sqrt {M}(1-M^{-1/4})^{-2}
			+\frac{4M^2}{M-1}c_1^4 c_3\right]\varepsilon^{-4}.$$
}

\section{Numerical Simulations}
To illustrate the theoretical results, we consider a non-linear scalar SDE
\begin{equation}\label{SDE2}
	dx(t)=(x(t)-x^3(t))dt+|x(t)|^{3/2}dB(t),\ \ \ t\geq 0,\ x(0)=x_0\in \mathbb{R}
\end{equation}
where $B(t)$ is a scalar Brownian motion. This is a specified Lewis stochastic volatility model. According to Examples 3.5 and 3.9 in \cite{Mao2016}, we sample over 1000
discretized Brownian paths and use stepsizes $s_l=T/2^l$ for $l=1,2,\dots,5$ in the truncated EM method. Let $\hat{Y_l}$ denote the sample value of $Y_l$.
Here we set $T=1$ and $h(s_l)=s_l^{-1/4}$.
\par
Firstly, we show some computational results of the classic EM method with the MLMC method.
\par
\begin{table}[!h]
\caption{Numerical results using the MLMC  with the classic EM method}
\end{table}
\begin{tabular}{|c|c|c|c|c|c|}
\hline
$l $ & 1 & 2 & 3 & 4 & 5 \\ \hline
$\hat{Y_l}$  & 1.00 & 2.59e+102 & -2.94e+159 & --- & --- \\ \hline
\end{tabular}
\medskip
\par \noindent
It can be seen from Table 1 that the simulation result of (\ref{SDE2}) computed by the MLMC approach together with the classic EM method is divergent.
\par
The simulation results using the MLMC method combined with the truncated EM method is presented in Table 2. It is clear that some convergent trend is displayed.
\par
\begin{table}[!h]
\caption{Numerical results using the MLMC  with the truncated EM method}
\end{table}
\begin{tabular}{|c|c|c|c|c|c|}
\hline
$l$ & 1 & 2 & 3 & 4 & 5 \\ \hline
$\hat{Y_l}$ & 0.39 & -0.18 & -0.024 & -0.003 & -0.0006 \\ \hline
\end{tabular}
\medskip
\par \noindent 
Next, it is noted that compared with the standard Monte Carlo method the computational cost can be saved by using MLMC method. From Figure 1, we can see that the MLMC method is approximately 10 times more efficient than the standard Monte Carlo method when $\varepsilon$ is sufficient small.
\begin{figure}[!h]
  \centering
  \includegraphics[scale=0.7]{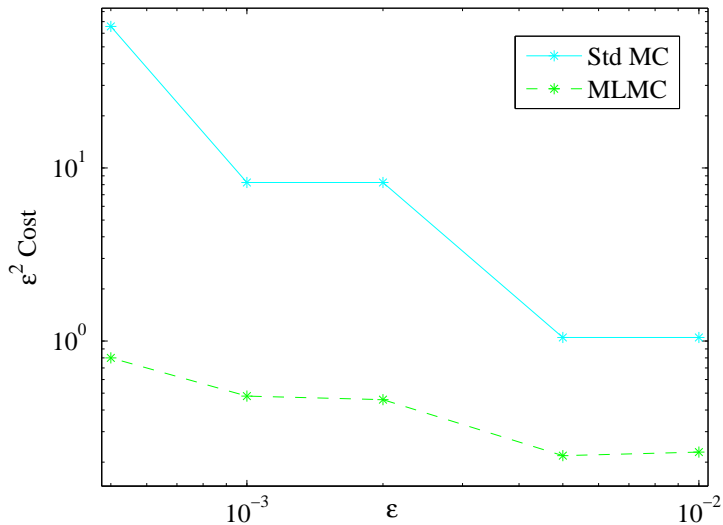}
  \caption{Computational cost}
\end{figure}

\section*{Acknowledgments}
The authors would like to thank the referee and editor for their very useful comments and suggestions, which have helped to improve this paper a lot.
\par
The authors would also like to thank
the Natural Science Foundation of Shanghai (14ZR1431300), 
the Royal Society (Wolfson Research Merit Award WM160014),
the Leverhulme Trust (RF-2015-385),
the Royal Society of London (IE131408),
the EPSRC (EP/E009409/1),
the Royal Society of Edinburgh (RKES115071),
the London Mathematical Society (11219),
the Edinburgh Mathematical Society (RKES130172),
the Ministry of Education (MOE) of China (MS2014DHDX020),
Shanghai Pujiang Program (16PJ1408000),
the Natural Science Fund of Shanghai Normal University (SK201603),
Young Scholar Training Program of Shanghai's Universities
for their financial support.

\section*{Appendix}
{
\it{Proof of Theorem \ref{generalthmMLMC}.}}

	Using the notation $\lceil x\rceil $ to denote the unique integer $n$ satisfying the inequalities $x\leq n<x+1$, we start by choosing $L$ to be
	\begin{equation*}
		L={\lceil \frac{log(\sqrt{2}c_1T^\alpha \varepsilon ^{-1})}{\alpha log M}\rceil},
	\end{equation*}
	so that
	\begin{equation*}
		\frac{1}{\sqrt{2}}M^{-\alpha}\varepsilon<c_1s_L^\alpha\leq \frac{1}{\sqrt{2}}\varepsilon.
	\end{equation*}
Hence, by the condition 1 and 2 we have
	\begin{equation}
		\begin{aligned}
			&\big(\mathbb{E}[Y]-\mathbb{E}[f(X(T))]\big)^2\\
=&\Big({\mathbb{E}}[\sum\limits_{l=0}^L Y_l]-{\mathbb{E}}[f(X(T))]\Big)^2\\
=& (\mathbb{E}[f(\overline{X}_{s_L}(T))-f(X(T))])^2 \\
\leq & (c_1s_L^\alpha)^2 \triangleq \frac{1}{2}\varepsilon^2.
		\end{aligned}
	\end{equation}
Therefore, we have
	$$(\mathbb{E}[Y]-\mathbb{E}[f(X)])^2\leq \frac{1}{2}\varepsilon ^2.$$
This upper bound on the square of bias error together with the upper bound of $\frac{1}{2}\varepsilon^2$ on the variance of the estimator, which will be proved later, gives a upper bound of $\varepsilon^2$ to the MSE.\par
Noting
	\begin{equation*}
		\sum\limits_{l=0}^Ls_l^{-1}=s_L^{-1}\sum\limits_{i=0}^LM^{-i}<\frac{M}{M-1}s_L^{-1},
	\end{equation*}
using the standard result for a geometric series and the inequality $\frac{1}{\sqrt{2}}M^{-\alpha}\varepsilon<c_1s_L^\alpha$, we can obtain
	\begin{equation*}
		s_L^{-1}<M(\frac{\varepsilon}{\sqrt{2}c_1})^{-1/{\alpha}}.
	\end{equation*}
Then, we have
	\begin{equation}
		\sum\limits_{l=0}^Ls_l^{-1}<\frac{M}{M-1}s_L^{-1}<\frac{M^2}{M-1}(\sqrt{2}c_1)^{\color{blue}{1/\alpha}}\varepsilon^{-1/\alpha}.
	\end{equation}\par
We now consider the different possible values of $\beta$ and to compare them to the $\alpha$.\par \noindent
	(a) If $\beta=1$, we set $N_l=\lceil 2\varepsilon^{-2}(L+1)c_2s_l\rceil$ so that
	\begin{equation*}
		V[Y]=\sum\limits_{l=0}^LV[Y_l]\leq\sum\limits_{l=0}^Lc_2N_l^{-1}s_l\leq\frac{1}{2}\varepsilon^2,
	\end{equation*}
which is the required.\par
For the bound of the computational complexity $C$, we have
	\begin{equation*}
		\begin{aligned}
			C&=\sum\limits_{l=0}^LC_l\leq c_3\sum\limits_{l=0}^LN_ls_l^{-1}\\
			&\leq c_3\sum\limits_{l=0}^L(2\varepsilon^{-2}(L+1)c_2s_l+1)s_l^{-1}\\
			&\leq c_3(2\varepsilon^{-2}(L+1)^2c_2+{\color{blue}\sum\limits_{l=0}^Ls_l^{-1}})\\
			&\leq c_3(2\varepsilon^{-2}(L+1)^2c_2+\frac{M^2}{M-1}(\sqrt{2}c_1)^{\color{blue}{1/\alpha}}\varepsilon^{-1/\alpha}).
		\end{aligned}
	\end{equation*}
According to the definition of $L$, we have
	\begin{equation*}
		L \leq {\frac{log\varepsilon^{-1}}{\alpha logM}+\frac{log(\sqrt{2}c_1T^\alpha)}{\alpha logM}+1}.
	\end{equation*}
	Given that $1<log\varepsilon^{-1}$ for $\varepsilon<e^{-1}$, we have
	\begin{equation*}
		L+1\leq c_5log\varepsilon^{-1},
	\end{equation*}
	where
	\begin{equation*}
		c_5=\frac{1}{\alpha logM}+max(0,\frac{log(\sqrt{2}c_1T^\alpha)}{\alpha logM})+2.
	\end{equation*}
	Hence, the computation complexity is bounded by
	\begin{equation*}
		\begin{aligned}
			C&\leq c_3(2\varepsilon^{-2}c_5^2(log\varepsilon^{-1})^2c_2+\frac{M^2}{M-1}(\sqrt{2}c_1)^{1/\alpha}\varepsilon^{-1/\alpha})\\
			&=c_3(2\varepsilon^{-2}c_5^2(log\varepsilon)^2c_2+\frac{M^2}{M-1}(\sqrt{2}c_1)^{1/\alpha}\varepsilon^{-1/\alpha}).
		\end{aligned}
	\end{equation*}
	So if $\alpha\leq (-log{\varepsilon})/{log[(log\varepsilon/\varepsilon)^2]}$, we have
	\begin{equation*}
		C\leq c_3(2c_5^2c_2+\frac{M^2}{M-1}(\sqrt{2}c_1)^{1/\alpha}){\varepsilon^{-1/\alpha}}.
	\end{equation*}
	If $\alpha> (-log{\varepsilon})/{log[(log\varepsilon/\varepsilon)^2]}$, we have
	\begin{equation*}
		C\leq c_3(2c_5^2c_2+\frac{M^2}{M-1}(\sqrt{2}c_1)^{1/\alpha}){\varepsilon^{-2}(log\varepsilon)^2}.
	\end{equation*}\par
\par \noindent
	(b) For $\beta>1$, setting
	\begin{equation*}
		N_l=\lceil 2\varepsilon ^{-2}c_2T^{(\beta-1)/2}(1-M^{-(\beta-1)/2})^{-1}s_l^{(\beta+1)/2}\rceil,
	\end{equation*}
	then we have
	\begin{equation*}
		\begin{aligned}
			V[Y]&=\sum\limits_{l=0}^{L}V[Y_l]\leq \sum\limits_{l=0}^{L}c_2N_l^{-1}s_l^\beta\\
			&\leq \frac{1}{2} \varepsilon^2T^{-(\beta-1)/2}(1-M^{-(\beta-1)/2})\sum\limits_{l=0}^Ls_l^{(\beta-1)/2}.
		\end{aligned}
	\end{equation*}
Using the stand result for a geometric series
	\begin{equation}
		\begin{aligned}
			\sum\limits_{l=0}^Ls_l^{(\beta-1)/2}&=T^{(\beta-1)/2}\sum\limits_{l=0}^L(M^{-(\beta-1)/2})^l\\
			&<T^{(\beta-1)/2}(1-M^{-(\beta-1)/2})^{-1},
		\end{aligned}
	\end{equation}
we obtain  that the upper bound of variance is $\frac{1}{2}\varepsilon^2$. So the computation complexity is bounded by
	\begin{equation*}
		\begin{aligned}
			C&\leq c_3\sum\limits_{l=0}^L N_ls_l^{-1}\\
			&\leq c_3\sum\limits_{l=0}^L (2\varepsilon^{-2}c_2T^{(\beta-1)/2}(1-M^{-(\beta-1)/2})^{-1} s_l^{(\beta+1)/2}+1)s_l^{-1}\\
			&=c_3[2\varepsilon^{-2}c_2T^{(\beta-1)/2}(1-M^{-(\beta-1)/2})^{-1}\sum\limits_{l=0}^L s_l^{(\beta-1)/2}+\sum\limits_{l=0}^L s_l^{-1}]\\
			&\leq c_3[2\varepsilon^{-2}c_2T^{(\beta-1)/2}(1-M^{-(\beta-1)/2})^{-1} T^{(\beta-1)/2}(1-M^{-(\beta-1)/2})^{-1}\\
			&+\frac{M^2}{M-1}(\sqrt{2}c_1)^{1/\alpha}\varepsilon^{-1/\alpha}]\\
			&=c_3[2\varepsilon^{-2}c_2T^{\beta-1}(1-M^{-(\beta-1)/2})^{-2}+\frac{M^2}{M-1}(\sqrt{2}c_1)^{1/\alpha}\varepsilon^{-1/\alpha}].
		\end{aligned}
	\end{equation*}
	So when $\alpha\geq\frac{1}{2}$, we have
	\begin{equation*}
		C\leq c_3[2c_2T^{\beta-1}(1-M^{-(\beta-1)/2})^{-2}+\frac{M^2}{M-1}(\sqrt{2}c_1)^{1/\alpha}]\varepsilon^{-2},
	\end{equation*}
When {\color{blue}$\alpha<\frac{1}{2}$,} we have
	\begin{equation*}
		C\leq c_3[2c_2T^{\beta-1}(1-M^{-(\beta-1)/2})^{-2}+\frac{M^2}{M-1}(\sqrt{2}c_1)^{1/\alpha}]\varepsilon^{-1/\alpha}.
	\end{equation*}
\par \noindent
{\color{blue}
	(c) For $0<\beta<1$, setting
	\begin{equation*}
		N_l=\lceil 2\varepsilon ^{-2}c_2s_L^{-(1-\beta)/2}(1-M^{-(1-\beta)/2})^{-1}s_l^{(\beta+1)/2}\rceil,
	\end{equation*}
	then we have
	\begin{equation*}
		\begin{aligned}
			V[Y]&=\sum\limits_{l=0}^{L}V[Y_l]\leq \sum\limits_{l=0}^{L}c_2N_l^{-1}s_l^\beta\\
			&\leq \frac{1}{2}\varepsilon ^{2}s_L^{(1-\beta)/2}(1-M^{-(1-\beta)/2})\sum\limits_{l=0}^{L}s_l^{-(1-\beta)/2}.
		\end{aligned}
	\end{equation*}
Because
\begin{equation}\label{62}
		\begin{aligned}
\sum\limits_{l=0}^{L}s_l^{-(1-\beta)/2}
&=s_L^{-(1-\beta)/2}\sum\limits_{l=0}^L(M^{-(1-\beta)/2})^l\\
&<s_L^{-(1-\beta)/2}(1-M^{-(1-\beta)/2})^{-1},
\end{aligned}
	\end{equation}
we obtain the upper bound on the variance of the estimator to be $\frac{1}{2}\varepsilon^2$.\par
Finally, using the upper bound of $N_l$ , the computational complexity is
	\begin{equation*}
		\begin{aligned}
			C&\leq c_3\sum\limits_{l=0}^L N_ls_l^{-1}\\
			&\leq c_3\sum\limits_{l=0}^L (2\varepsilon^{-2}c_2s_L^{-(1-\beta)/2}(1-M^{-(1-\beta)/2})^{-1} s_l^{(\beta+1)/2}+1)s_l^{-1}\\
			&=c_3[2\varepsilon^{-2}c_2s_L^{-(1-\beta)/2}(1-M^{-(1-\beta)/2})^{-1}\sum\limits_{l=0}^L s_l^{-(1-\beta)/2}+\sum\limits_{l=0}^L s_l^{-1}]\\
			&\leq c_3[2\varepsilon^{-2}c_2s_L^{-(1-\beta)}(1-M^{-(1-\beta)/2})^{-2} +\sum\limits_{l=0}^L s_l^{-1}],
		\end{aligned}
	\end{equation*}
where (\ref{62}) is used in the last inequality.\par
Moreover, because of the inequality $\frac{1}{\sqrt{2}}M^{-\alpha}\varepsilon<c_1s_L^\alpha$, we have
	\begin{equation*}
		s_L^{-(1-\beta)}<(\sqrt{2}c_1)^{(1-\beta)/\alpha}M^{1-\beta}\varepsilon^{-(1-\beta)/\alpha},
	\end{equation*}
}
	then
	\begin{equation*}
		\begin{aligned}
			C&\leq c_3[2\varepsilon^{-2}c_2s_L^{-(1-\beta)}(1-M^{-(1-\beta)/2})^{-2} +\sum\limits_{l=0}^L s_l^{-1}]\\
			&\leq c_3[2\varepsilon^{-2}c_2(\sqrt{2}c_1)^{(1-\beta)/\alpha}M^{1-\beta}\varepsilon^{-(1-\beta)/\alpha}(1-M^{-(1-\beta)/2})^{-2}
			+\sum\limits_{l=0}^L s_l^{-1}]\\
			&\leq c_3[2\varepsilon^{-2}c_2(\sqrt{2}c_1)^{(1-\beta)/\alpha}M^{1-\beta}\varepsilon^{-(1-\beta)/\alpha}(1-M^{-(1-\beta)/2})^{-2}
			+\frac{M^2}{M-1}(\sqrt{2}c_1)^{1/\alpha}\varepsilon^{-1/\alpha}]\\
			&=c_3[2c_2(\sqrt{2}c_1)^{(1-\beta)/\alpha}M^{1-\beta}(1-M^{-(1-\beta)/2})^{-2}\varepsilon^{-2-(1-\beta)/\alpha}
			+\frac{M^2}{M-1}(\sqrt{2}c_1)^{1/\alpha}\varepsilon^{-1/\alpha}].
		\end{aligned}
	\end{equation*}
	If $\beta\leq 2\alpha$, then
	$\varepsilon^{-2-(1-\beta)/\alpha}>\varepsilon^{-\frac{1}{\alpha}}$,
	so we have
	\begin{equation*}
		C\leq c_3[2c_2(\sqrt{2}c_1)^{(1-\beta)/\alpha}M^{1-\beta}(1-M^{-(1-\beta)/2})^{-2}
		+\frac{M^2}{M-1}(\sqrt{2}c_1)^{1/\alpha}]\varepsilon^{-2-(1-\beta)/\alpha}.
	\end{equation*}
	If $\beta> 2\alpha$, then
	$\varepsilon^{-2-(1-\beta)/\alpha}<\varepsilon^{-\frac{1}{\alpha}}$, so we have
	\begin{equation*}
		C\leq c_3[2c_2(\sqrt{2}c_1)^{(1-\beta)/\alpha}M^{1-\beta}(1-M^{-(1-\beta)/2})^{-2}
		+\frac{M^2}{M-1}(\sqrt{2}c_1)^{1/\alpha}]\varepsilon^{-1/\alpha}.
	\end{equation*}

\bibliographystyle{plain}


\end{document}